\documentclass[10pt,a4paper]{amsart}

\usepackage{pst-plot}
   \usepackage{pstricks}
   \usepackage{multido}
   \usepackage{pstricks-add} 

\usepackage{url}    
\usepackage{amsmath}
\usepackage{amsfonts}
\usepackage{color} 

\newtheorem{theo}{Theorem}[section]

\begin{document}

\date{\today}

\title[Abelian  varieties]{The characteristic polynomials of abelian varieties of dimensions 4 over finite fields}
\author{Safia Haloui, Vijaykumar Singh}
\address{Institut de Math\'ematiques de Luminy, Marseille, France\\
and\\
Claude Shannon Institute, Dublin, Irland}

\email{haloui@iml.univ-mrs.fr,vijaykumar.singh@ucdconnect.ie}



\begin{abstract}
We describe the set  of characteristic polynomials of abelian varieties of dimension $4$ over finite fields.
\end{abstract}

\maketitle
\section{Introduction and results}
The aim of this paper is to give a description of the set of characteristic polynomials of abelian varieties of dimension $4$ as it was done in \cite{saf} for dimension $3$.

\bigskip

It is well known that the characteristic polynomial of an abelian variety of dimension $g$ over $\mathbb{F}_q$ (with $q=p^n$) is monic, with integer coefficients, of degree $2g$ and that the sets of its roots consist of couples of complex conjugated numbers of modulus $\sqrt{q}$.  Any polynomial having those properties is called a \textit{Weil polynomial}.  Obviously, every Weil polynomial of degree $8$ has the form
$$p(t)=t^8+a_1t^7+a_2t^6+a_3t^5+a_4t^4+qa_3t^3+q^2a_2t^2+q^3a_1t+q^4$$
for certain integers $a_1$, $a_2$, $a_3$ and $a_4$. In  Section 2 we prove the following proposition which gives a characterization of the quadruples $(a_1 , a_2 , a_3, a_4 )$ corresponding to Weil polynomials of degree $8$ (see \cite{water, ruck, saf}  for a characterization of Weil polynomials of lower degrees):
\begin{theo}\label{coeffweil}
Let $p(t)=t^8+a_1t^7+a_2t^6+a_3t^5+a_4t^4+qa_3t^3+q^2a_2t^2+q^3a_1t+q^4$ be a polynomial with integer coefficients. Then $p(t)$ is a Weil polynomial if and only if
  either
$$p(t)=(t^2\pm\sqrt{q})^2h(t)$$
where $h(t)$ is a Weil polynomial, or
 the following conditions hold:
\begin{enumerate}
\item $\vert a_1\vert< 8\sqrt q$,
\item $6\sqrt q\vert a_1\vert -20q< a_2\leq\frac{3a_1^2}{8}+4q$,
\item $-9qa_1-4\sqrt qa_2-16q\sqrt q< a_3< -9qa_1+4\sqrt qa_2+16q\sqrt q$,
\item $-\frac{a_1^3}{8}+\frac{a_1a_2}{2}+qa_1-(\frac{2}{3}(\frac{3a_1^2}{8}-a_2+4q))^{3/2}\leq a_3\leq -\frac{a_1^3}{8}+\frac{a_1a_2}{2}+qa_1+(\frac{2}{3}(\frac{3a_1^2}{8}-a_2+4q))^{3/2}$,
\item $2\sqrt q\vert qa_1+a_3\vert-2qa_2-2q^2< a_4$,
\item $\frac{9a_1^4}{256}-\frac{3a_1^2a_2}{16}+\frac{a_1a_3}{4}+\frac{a_2^2}{6}+\frac{2qa_2}{3}+\frac{2q^2}{3}+\omega+\overline{\omega}\leq a_4\leq\frac{9a_1^4}{256}-\frac{3a_1^2a_2}{16}+\frac{a_1a_3}{4}+\frac{a_2^2}{6}+\frac{2qa_2}{3}+\frac{2q^2}{3}+j\omega+j^2\overline{\omega}$
\end{enumerate}
where \\
$\omega=\frac{1}{24}\left (8(-\frac{3 a_1^2}{8}+ a_2-4q)^6+540(-\frac{3 a_1^2}{8}+ a_2-4q)^3(\frac{a_1^3}{8}-qa_1-\frac{a_1a_2}{2}+a_3)^2 \right . $ \\ 
$\left . -729(\frac{a_1^3}{8}-qa_1-\frac{a_1a_2}{2}+a_3)^4\right . $ \\ 
${\left . +i9\vert \frac{a_1^3}{8}-qa_1-\frac{a_1a_2}{2}+a_3\vert(-(\frac{a_1^3}{8}-qa_1-\frac{a_1a_2}{2}+a_3 )^2-\frac{8}{27}(-\frac{3 a_1^2}{8}+ a_2-4q)^3)^{3/2}\right )}^{1/3},$\\
$\omega^{1/3} = \vert\omega\vert^{1/3}e^{\frac{arg(\omega)i}{3}}$ and $j=e^{\frac{2i\pi}{3}}$.
\end{theo}

Now it remains to give criterions to determine when a Weil polynomial is the characteristic polynomial of an abelian variety of dimension $4$.  Since the characteristic polynomial of a non-simple abelian variety is a product of characteristic polynomials of abelian varieties of smaller dimensions and our problem is already solved for smaller dimensions \cite{water, ruck, saf}, it is sufficient to consider the simple case.

By results of Honda and Tate, the characteristic polynomial of a simple abelian variety of dimension $4$ over $\mathbb{F}_q$ has the form $p(t)=h(t)^e$ where $h(t)$ is an irreducible  Weil polynomial and $e$ is an integer.  Xing \cite{xing} and Maisner and Nart \cite{mana} gave independently a description of characteristic polynomials of abelian varieties of dimension $4$ with $e>1$. Therefore, we can restrict our attention to the case $e=1$, that is, $p(t)$ is irreducible.

\bigskip

If $p(t)$ is irreducible, the determination of the possible Newton polygons for $p(t)$ (Section 3) gives us the following proposition:
\begin{theo}\label{newt} Let $p(t)=t^8+a_1t^7+a_2t^6+a_3t^5+a_4t^4+qa_3t^3+q^2a_2t^2+q^3a_1t+q^4$ be an irreducible Weil polynomial. Then $p(t)$ is the characteristic polynomial of an abelian variety of dimension $4$ if and only if one of the following conditions holds:
\begin{enumerate}
\item $v_p(a_4)=0$,
\item $v_p(a_3)=0$, $v_p(a_4)\geq n/2$ and $p(t)$ has no root of valuation $n/2$ in $\mathbb{Q}_p$,
\item $v_p(a_2)=0$, $v_p(a_3)\geq n/2$, $v_p(a_4)\geq n$ and $p(t)$ has no root of valuation $n/2$ in $\mathbb{Q}_p$,
\item $v_p(a_1)=0$, $v_p(a_2)\geq n/2$, $v_p(a_3)\geq n$, $v_p(a_4)\geq 2n$ and $p(t)$ has no root of valuation $n/2$ nor factor of degree $3$ in $\mathbb{Q}_p$,
\item $v_p(a_1)=0$, $v_p(a_2)\geq n/3$, $v_p(a_3)\geq 2n/3$, $v_p(a_4)=n$ and $p(t)$ has no root of valuation $n/3$ and $2n/3$ in $\mathbb{Q}_p$,
\item $v_p(a_1)\geq n/3$, $v_p(a_2)\geq 2n/3$, $v_p(a_3)=n$, $v_p(a_4)\geq 3n/2$ and $p(t)$ has no root in $\mathbb{Q}_p$,
\item $v_p(a_1)\geq n/4$, $v_p(a_2)\geq n/2$, $v_p(a_3)=3n/4$, $v_p(a_4)=n$ and $p(t)$ has no root nor factor of degree $2$ and $3$ in $\mathbb{Q}_p$,
\item $v_p(a_1)\geq n/2$, $v_p(a_2)\geq n$, $v_p(a_3)=3n/2$, $v_p(a_4)\geq 2n$ and $p(t)$ has no root  nor factor of degree $3$ in $\mathbb{Q}_p$.
\end{enumerate}
The $p$-ranks of abelian varieties in cases (1), (2), (3), (4), (5), (6), (7) and (8) are respectively  $4$, $3$, $2$, $1$, $1$, $0$, $0$ and $0$. The abelian varieties in case (8) are supersingular.
\end{theo}

It is possible to make condition (8) of Theorem \ref{newt} more explicit. 
Indeed, in \cite{Vijay}, Singh, McGuire and Zaytsev gave the list of irreducible characteristic polynomials of  supersingular abelian varieties of dimension $4$ , where $q$ is not a square. We complete the classification by finding the list in the case $q$ is a square in the following proposition (see Section 4).

\begin{theo}The polynomial $p(t)$ is the irreducible characteristic polynomial of a supersingular abelian
variety of dimension $4$  if and only if one of the following conditions holds
\begin{itemize}
\item $q$ is a square and $(a_1 , a_2 , a_3, a_4 )$ belongs to the following list:
\begin{enumerate}
\item  $( -q^{1/2},0 , q^{3/2},-q^2 ), ~p \not\equiv 1 \mod 15 $,
\item  $( q^{1/2},0 , -q^{3/2},-q^2 ),~ p \not\equiv 1 \mod 30 $,
 \item $(0,0,0,0),~ p ~\not\equiv 1 \mod 16 $,
\item $(0,-q,0,q^2),~ p \not\equiv 1 \mod 20$,  
\item $(0,0,0,-q^2), ~p \not \equiv 1 \mod 24 $,
\end{enumerate}
\item $q$ is not a square  and $(a_1 , a_2 , a_3, a_4 )$ belongs to the following list:

\begin{enumerate}
\item $(\pm \sqrt{pq}, q, 0, -q^2), ~p =2$,
\item $(\pm \sqrt{pq},2q,\pm q\sqrt{pq},q^2 ), ~p=3$,
\item $(0,0,0,0)$,
\item $(0,-q, 0, q^2)$,
\item $(0,q, 0, q^2),~ p \neq 5$,
\item $(0,0,0,-q^2),~ p \neq 2$,
\item $(\pm \sqrt{pq},2q,\pm q\sqrt{pq},3q^2 ),~ p=5$.
\end{enumerate}
\end{itemize}
\end{theo}


\section{The coefficients of Weil polynomials of degree $8$}
In this section, we prove Theorem \ref{coeffweil}. It is clear that a Weil polynomial with a real root must have the form 
$$p(t)=(t^2\pm\sqrt{q})^2h(t)$$
where $h(t)$ is a Weil polynomial. Conversely, these polynomials are Weil polynomials.

\bigskip

Let $p(t)=t^8+a_1t^7+a_2t^6+a_3t^5+a_4t^4+qa_3t^3+q^2a_2t^2+q^3a_1t+q^4\in\mathbb{Z}[t]$ be a polynomial with no real root. Then the set of the roots of $p(t)$ consists of pairs of complex conjugated numbers, say $\omega_1,\overline{\omega_1},\dots ,\omega_4, \overline{\omega_4}$. Letting $x_i=-(\omega_i+\overline{\omega_i})$ we have $p(t)=\prod_{i=1}^4(t^2+x_it+q)$. Arguing as in \cite{saf}, $p(t)$ is a Weil polynomial if and only if  the polynomials $f^+(t)=\prod_{i=1}^4(t-(2\sqrt{q}+x_i))$ and $f^-(t)=\prod_{i=1}^4(t-(2\sqrt{q}-x_i))$ have only real and positive roots.

\bigskip

First, we determine a necessary and sufficient condition of some polynomial of degree $4$ having all real roots. 

Let $f(t)=t^4+r_1t^3+r_2t^2+r_3t+r_4$ be a monic polynomial of degree $4$ with real coefficients.
Looking at the table of variation of $f(t)$, we see that there exists some $r_4$ for which  $f(t)$ has all real roots if and only if $f'(t)$ has all real roots. This condition is equivalent to
\begin{eqnarray}\label{discder}
\Delta_{f'}\geq 0
\end{eqnarray}
where $\Delta_{f'}$ is the discriminant of $f'(t)$.

The discriminant of $f(t)$ is a polynomial of degree $3$ in $r_4$ which we will denote $\Delta_{f}(t)$ (that is, $\Delta_{f} (r_4)$ is the discriminant of $f(t)$). It is well known that if  $f(t)$ has all real roots  then $\Delta_{f} (r_4)\geq 0$. Moreover, the function which associate to $r_4$ the number of roots of $f(t)$ is constant  on the intervals delimited by the roots of $\Delta_{f} (t)$ (because $\Delta_{f} (r_4)=0$ when $f(t)$ has a multiple root). 

When $r_4$ is very big, the graph of $f(t)$ doesn't touch the $x$-axis and therefore $f(t)$ has no real roots. Thus if $\gamma_3$ is the biggest root of $\Delta_{f} (t)$, by the previous discussion, $f(t)$ has no real root for $r_4\in ]\gamma_3 ;+\infty [$.

We deduce that if (\ref{discder}) is satisfied then $\Delta_{f}(t)$ must have $3$ real roots $\gamma_1\leq\gamma_2\leq\gamma_3$ and $f(t)$ has all real roots if and only if 
\begin{eqnarray}\label{r4}
\gamma_1\leq  r_4\leq\gamma_2.
\end{eqnarray}

The roots of $\Delta_{f}(t)$ can be found using Cardan's method. Let us recall quickly what it is. 

Given a polynomial $h(t)=t^3+u_2t+u_3$, we set $\delta=-u_3^2-\frac{4}{27}u_2^3$. Then $h(t)$ has all real roots if and only if $\delta \geq 0$. If this is the case, the roots of $h(t)$ are $\gamma_1=\omega+\overline{\omega}$, $\gamma_2=j\omega+j^2\overline{\omega}$ and $\gamma_3=j^2\omega+j\overline{\omega}$ where $j=e^{\frac{2i\pi}{3}}$ and $\omega= ( \frac{-u_3+i\sqrt{\delta}}{2} ) ^{1/3}$. Moreover, with the convention $\omega^{1/3} = \vert\omega\vert^{1/3}e^{\frac{arg(\omega)i}{3}}$, we have $\gamma_1\leq\gamma_2\leq\gamma_3$.

In the general case, we have $h(t)=v_0t^3+v_1t^2+v_2t+v_3$ and we conclude using the fact that $\frac{1}{v_0}h(t-\frac{v_1}{3v_0})=t^3+u_2t+u_3$ with $u_2=-\frac{v_1^2}{3v_0^2}+\frac{v_2}{v_0}$ and $u_3=\frac{2v_1^3}{27v_0^3}-\frac{v_1v_2}{3v_0^2}+\frac{v_3}{v_0}$.

\bigskip

For $i=1,2,3,4$, let  $s_i$ denote the $i$th symmetric function of the $x_i$'s (that is, $\prod_{i=1}^{4}(t+x_i)=t^4+\sum_{i=1}^{4}s_it^{4-i})$. Expanding the expression $p(t)=\prod_{i=1}^4(t^2+x_it+q)$, we find:
\begin{eqnarray*}
s_1 & = & a_1\\
s_2 & = & a_2-4q\\
s_3 & = & a_3-3qa_1\\
s_4 & = & a_4-2qa_2+2q^2.
\end{eqnarray*}

Now, in order to simplify the calculation, we remark that $f^+(t)$ and $f^-(t)$  have all real roots if and only if the polynomial $f(t)=\prod_{i=1}^4(t+x_i-\frac{a_1}{4})$ has. Therefore, it is equivalent to apply our results to $f(t)$. 

Expanding the expression of $f(t)$,  we find that $f(t)=t^4+r_2t^2+r_3t+r_4$, where
\begin{eqnarray*}
r_2 & = & -\frac{3 s_1^2}{8}+ s_2\\
r_3 & = & \frac{s_1^3}{8}-\frac{s_1s_2}{2}+s_3\\
r_4 & = & -\frac{3 s_1^4}{256}+\frac{s_1^2 s_2}{16}-\frac{s_1s_3}{4}+ s_4.
\end{eqnarray*}

Substituting $s_1$, $s_2$, $s_3$ and $s_4$ with their expressions in $a_1$, $a_2$, $a_3$ and $a_4$ we obtain
\begin{eqnarray*}
r_2 & = & -\frac{3 a_1^2}{8}+ a_2-4q\\
r_3 & = & \frac{a_1^3}{8}-qa_1-\frac{a_1a_2}{2}+a_3\\
r_4 & = & -\frac{3 a_1^4}{256}+\frac{qa_1^2 }{2}+\frac{a_1^2 a_2}{16}-\frac{a_1a_3}{4} -2qa_2+2q^2+ a_4.
\end{eqnarray*}

\bigskip

We have
$$\Delta_{f}(t)=256t^3-128r_2^2t^2+16r_2(r_2^3+9r_3^2)t-r_3^2(4r_2^3+27r_3^2).$$
Now, we use Cardan's method. Set
\begin{eqnarray*}
u_2 & = & -\frac{r_2^4}{48} + \frac{9 r_2r_3^2}{16}\\
u_3 & = & \frac{r_2^6}{864} + \frac{5 r_2^3 r_3^2}{64} - \frac{27r_3^4}{256}\\
\delta & = & -u_3^2-\frac{4}{27}u_2^3\; =\;\frac{r_3^2 (-8 r_2^3 - 27 r_3^2)^3}{1769472}.
\end{eqnarray*}
Suppose that (\ref{discder}) is satisfied. Then $\delta\geq 0$ and the roots of $\Delta_{f}(t)$ are
\begin{eqnarray*}
\gamma_1 & = & \omega+\overline{\omega}+\frac{r_2^2}{6}\\
\gamma_2 & = & j\omega+j^2\overline{\omega}+\frac{r_2^2}{6}\\
\gamma_3 & = & j^2\omega+j\overline{\omega}+\frac{r_2^2}{6}.
\end{eqnarray*}
where $$\omega=\frac{1}{24}{\left ( 8r_2^6+540r_2^3r_3^2-729r_3^4+i9\vert r_3\vert(-r_3^2-\frac{8}{27}r_2^3)^{3/2}\right )}^{1/3}.$$
Substituting $r_2$, $r_3$ and $r_4$ with their expressions in $a_1$, $a_2$, $a_3$ and $a_4$ we obtain condition (6) of Theorem \ref{coeffweil}.

Next we have to determine when (\ref{discder}) is satisfied. We have: 
$$\Delta_{f'}=-16(8 r_2^3 + 27 r_3^2).$$
Therefore, (\ref{discder}) is equivalent to
$$r_2\leq 0\; \mbox{ and }\; -{(\frac{-2r_2}{3})}^{3/2}\leq r_3\leq{(\frac{-2r_2}{3})}^{3/2}.$$
This gives us the second inequality of condition (2) and condition (4) of Theorem \ref{coeffweil}.

\bigskip

Finally, we determine when  the polynomials $f^+(t)$ and $f^-(t)$ have only positive roots. 

For $i=1,2,3$, let  $r^+_i$ and $r^-_i$ denote the respective $i$th coefficients of $f^+(t)$ and $f^-(t)$. Expending the expressions $f^+(t)=\prod_{i=1}^4(t-(2\sqrt{q}+x_i))$ and $f^-(t)=\prod_{i=1}^4(t-(2\sqrt{q}-x_i))$, we find:
\begin{eqnarray*}
r^+_1 & = & -8\sqrt q-s_1\\
r^+_2 & = & 24q+6\sqrt qs_1+s_2\\
r^+_3 & = & -32q\sqrt q-12qs_1-4\sqrt qs_2-s_3\\
r^+_4 & = & 16q^2+8q\sqrt qs_1+4qs_2+2\sqrt qs_3+s_4
\end{eqnarray*}
and
\begin{eqnarray*}
r^-_1 & = & -8\sqrt q+s_1\\
r^-_2 & = & 24q-6\sqrt qs_1+s_2\\
r^-_3 & = & -32q\sqrt q+12qs_1-4\sqrt qs_2+s_3\\
r^-_4 & = & 16q^2-8q\sqrt qs_1+4qs_2-2\sqrt qs_3+s_4.
\end{eqnarray*}
Substituting $s_1$, $s_2$, $s_3$ and $s_4$ with their expressions in $a_1$, $a_2$, $a_3$ and $a_4$ we obtain
\begin{eqnarray*}
r^+_1 & = & -8\sqrt q-a_1 \\
r^+_2 & = & 20q+6\sqrt qa_1+a_2\\
r^+_3 & = & -16q\sqrt q-9qa_1-4\sqrt qa_2-a_3\\
r^+_4 & = & 2q^2+2q\sqrt qa_1+2qa_2+2\sqrt qa_3+a_4
\end{eqnarray*}
and
\begin{eqnarray*}
r^-_1 & = & -8\sqrt q+a_1\\
r^-_2 & = & 20q-6\sqrt qa_1+a_2 \\
r^-_3 & = & -16q\sqrt q+9qa_1-4\sqrt qa_2+a_3\\
r^-_4 & = & 2q^2-2q\sqrt qa_1+2qa_2-2\sqrt qa_3+a_4.
\end{eqnarray*}
Suppose that $f^+(t)$ and $f^-(t)$ have all real roots. Then by \cite[\S 2, Lemma]{smyth}, $f^+(t)$ and $f^-(t)$ have only positive roots if and only if $(-1)^ir^+_i>0$ and $(-1)^ir^-_i>0$ for $i=1,2,3,4$. This gives us the remaining conditions of  Theorem \ref{coeffweil} and concludes the proof.

\bigskip
\noindent
{\bf Remark.}  We could have used  \cite[\S 2, Lemma]{smyth} to determine when a polynomial of degree $4$ has only real roots but the computation and the results would have been longer.

\section{Newton polygons}
Let $p(t)$ be an irreducible Weil polynomial. By \cite{milwat}, $p(t)^e$ is the characteristic polynomial of a simple abelian variety, where $e$ the least common denominator of $v_p(f(0))/n$ where $f(t)$ runs through the irreducible factors of $p(t)$ over $\mathbb{Q}_p$. Thus $p(t)$ is the characteristic polynomial of an abelian variety of dimension $4$ if and only if  $e$ is equal to $1$ that is, $v_p(f(0))/n$ are integers.

 In order to determine when this condition is satisfied, we consider the Newton polygon of $p(t)$ (see \cite{weiss}).
Each of its edges define a factor of $p(t)$ over $\mathbb{Q}_p$. The degree of this factor is the length of the projection onto the $x$-axis of the edge and all the roots of this factor have a valuation equal to the slope of the edge. Therefore $e=1$ implies that the length of the projection onto the $x$-axis of any edge times its slope is a multiple of $n$. 

We graph the Newton polygons satisfying this condition and in each case, we give a necessary and sufficient condition to have $e=1$. The obtained results are summarized in Theorem \ref{newt}.

\bigskip
\noindent
\textbf{Ordinary case: }$\mathbf{v_p(a_4)=0}$

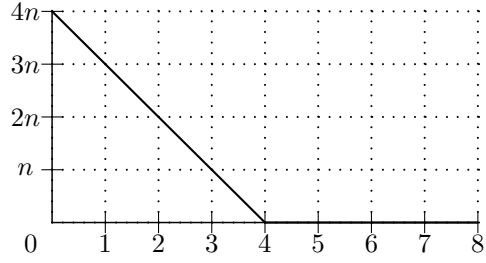
\begin{figure}[h]
   \centering
  \begin{pspicture}(0,0)(8,4)
  \psset{unit=7mm,labels=none}
  \rput(-0.4,-0.4){$0$}
  \rput(-0.5,1){$n$}
  \rput(-0.5,2){$2n$}
  \rput(-0.5,3){$3n$}
  \rput(-0.5,4){$4n$}
  \rput(1,-0.4){$1$}
    \rput(2,-0.4){$2$}
    \rput(3,-0.4){$3$}
   \rput(4,-0.4){$4$}
    \rput(5,-0.4){$5$}
    \rput(6,-0.4){$6$}
    \rput(7,-0.4){$7$}
    \rput(8,-0.4){$8$}
   \psaxes[linewidth=.7\pslinewidth]{-}(0,0)(0,0)(8,4)
  \psgrid[griddots=5, subgriddiv=0, gridlabels=0pt](0,0)(8,4)
   \psline(0,4)(4,0)(8,0)
\end{pspicture}
 \caption{\label{ordinaire}Ordinary case}
\end{figure}

The Newton polygon of $p(t)$ is represented in Figure \ref{ordinaire} and we always have $e=1$.

\bigskip
\noindent
\textbf{p-rank 3 case: }$\mathbf{v_p(a_4)>0}$\textbf{ and }$\mathbf{v_p(a_3)=0}$

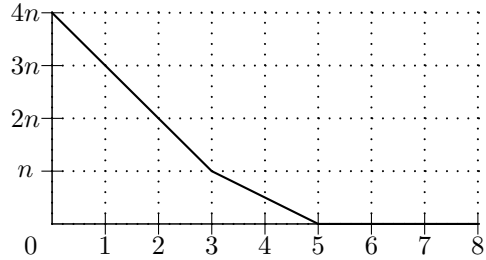
\begin{figure}[h]
   \centering
     \begin{pspicture}(0,0)(8,4)
  \psset{unit=7mm,labels=none}
    \rput(-0.4,-0.4){$0$}
  \rput(-0.5,1){$n$}
  \rput(-0.5,2){$2n$}
  \rput(-0.5,3){$3n$}
  \rput(-0.5,4){$4n$}
  \rput(1,-0.4){$1$}
    \rput(2,-0.4){$2$}
    \rput(3,-0.4){$3$}
   \rput(4,-0.4){$4$}
    \rput(5,-0.4){$5$}
    \rput(6,-0.4){$6$}
    \rput(7,-0.4){$7$}
    \rput(8,-0.4){$8$}
   \psaxes[linewidth=.7\pslinewidth]{-}(0,0)(0,0)(8,4)
  \psgrid[griddots=5, subgriddiv=0, gridlabels=0pt](0,0)(8,4)
   \psline(0,4)(3,1)(5,0)(8,0)
\end{pspicture}
 \caption{\label{prang3}$p$-rank $3$ case}
\end{figure}

The only Newton polygon for which $e=1$ is represented in Figure \ref{prang3}.

This is the Newton polygon of  $p(t)$ if and only if $v_p(a_4)\geq n/2$. If this condition holds, $p(t)$ has a factor  in $\mathbb{Q}_p$ of degree $2$ with roots of valuation $n/2$ and $e=1$ if and only if this factor is irreducible, that is, if and only if $p(t)$ has no root of valuation $n/2$ in $\mathbb{Q}_p$.


\bigskip
\noindent
\textbf{p-rank 2 case: }$\mathbf{v_p(a_4)>0}$\textbf{, }$\mathbf{v_p(a_3)>0}$\textbf{ and }$\mathbf{v_p(a_2)=0}$

\begin{figure}[h]
   \centering
    \begin{pspicture}(0,0)(8,4)
  \psset{unit=7mm,labels=none}
   \rput(-0.4,-0.4){$0$}
  \rput(-0.5,1){$n$}
  \rput(-0.5,2){$2n$}
  \rput(-0.5,3){$3n$}
  \rput(-0.5,4){$4n$}
  \rput(1,-0.4){$1$}
    \rput(2,-0.4){$2$}
    \rput(3,-0.4){$3$}
   \rput(4,-0.4){$4$}
    \rput(5,-0.4){$5$}
    \rput(6,-0.4){$6$}
    \rput(7,-0.4){$7$}
    \rput(8,-0.4){$8$}
   \psaxes[linewidth=.7\pslinewidth]{-}(0,0)(0,0)(8,4)
  \psgrid[griddots=5, subgriddiv=0, gridlabels=0pt](0,0)(8,4)
   \psline(0,4)(2,2)(4,1)(6,0)(8,0)
\end{pspicture}
 \caption{\label{prang2}$p$-rank $2$ case}
\end{figure}
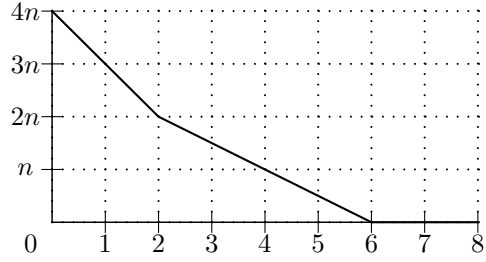

The only Newton polygon for which $e=1$ is represented in Figure \ref{prang2}.

This is the Newton polygon of  $p(t)$ if and only if $v_p(a_3)\geq n/2$ and $v_p(a_4)\geq n$. If these conditions hold, $p(t)$ has a factor  in $\mathbb{Q}_p$ of degree $4$ with roots of valuation $n/2$ and $e=1$ if and only if this factor has no root in $\mathbb{Q}_p$, that is, if and only if $p(t)$ has no root of valuation $n/2$ in $\mathbb{Q}_p$.

\bigskip
\noindent
\textbf{p-rank 1 case: }$\mathbf{v_p(a_4)>0}$\textbf{, }$\mathbf{v_p(a_3)>0}$\textbf{, }$\mathbf{v_p(a_2)>0}$\textbf{ and }$\mathbf{v_p(a_1)=0}$

\begin{figure}[h]
   \centering
  \begin{pspicture}(0,0)(8,4)
  \psset{unit=7mm,labels=none}
    \rput(-0.4,-0.4){$0$}
  \rput(-0.5,1){$n$}
  \rput(-0.5,2){$2n$}
  \rput(-0.5,3){$3n$}
  \rput(-0.5,4){$4n$}
  \rput(1,-0.4){$1$}
    \rput(2,-0.4){$2$}
    \rput(3,-0.4){$3$}
   \rput(4,-0.4){$4$}
    \rput(5,-0.4){$5$}
    \rput(6,-0.4){$6$}
    \rput(7,-0.4){$7$}
    \rput(8,-0.4){$8$}
   \psaxes[linewidth=.7\pslinewidth]{-}(0,0)(0,0)(8,4)
  \psgrid[griddots=5, subgriddiv=0, gridlabels=0pt](0,0)(8,4)
   \psline(0,4)(1,3)(7,0)(8,0)
\end{pspicture}
 \caption{\label{prang11}$p$-rank $1$ first case}
\end{figure}
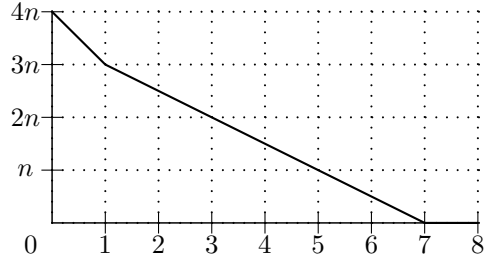

There are two Newton polygons for which $e=1$. One is represented in Figure \ref{prang11}.

This is the Newton polygon of  $p(t)$ if and only if $v_p(a_2)\geq n/2$, $v_p(a_3)\geq n$ and $v_p(a_4)\geq 2n$.  If these conditions hold,  $e=1$ if and only if $p(t)$ has no root of valuation $n/2$  nor factor of degree $3$ in $\mathbb{Q}_p$.

\medskip

\begin{figure}[h]
   \centering
    \begin{pspicture}(0,0)(8,4)
  \psset{unit=7mm,labels=none}
   \rput(-0.4,-0.4){$0$}
  \rput(-0.5,1){$n$}
  \rput(-0.5,2){$2n$}
  \rput(-0.5,3){$3n$}
  \rput(-0.5,4){$4n$}
  \rput(1,-0.4){$1$}
    \rput(2,-0.4){$2$}
    \rput(3,-0.4){$3$}
   \rput(4,-0.4){$4$}
    \rput(5,-0.4){$5$}
    \rput(6,-0.4){$6$}
    \rput(7,-0.4){$7$}
    \rput(8,-0.4){$8$}
   \psaxes[linewidth=.7\pslinewidth]{-}(0,0)(0,0)(8,4)
  \psgrid[griddots=5, subgriddiv=0, gridlabels=0pt](0,0)(8,4)
   \psline(0,4)(1,3)(4,1)(7,0)(8,0)
\end{pspicture}
 \caption{\label{prang12}$p$-rank $1$ second case}
\end{figure}
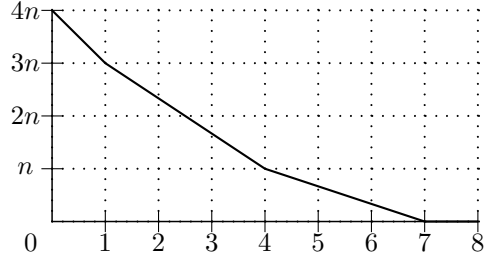

The other Newton polygon is represented in Figure \ref{prang12}.

This is the Newton polygon of  $p(t)$ if and only if $v_p(a_2)\geq n/3$, $v_p(a_3)\geq 2n/3$ and $v_p(a_4)=n$.  If these conditions hold $e=1$ if and only if $p(t)$ has no root of valuation $n/3$ and $2n/3$ in $\mathbb{Q}_p$.

\bigskip
\noindent
\textbf{p-rank 0 case: }$\mathbf{v_p(a_4)>0}$\textbf{, }$\mathbf{v_p(a_3)>0}$\textbf{, }$\mathbf{v_p(a_2)>0}$\textbf{ and }$\mathbf{v_p(a_1)>0}$

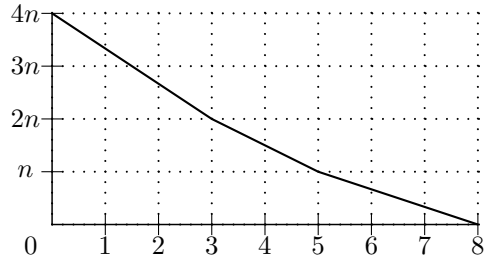
\begin{figure}[h]
   \centering
      \begin{pspicture}(0,0)(8,4)
  \psset{unit=7mm,labels=none}
    \rput(-0.4,-0.4){$0$}
  \rput(-0.5,1){$n$}
  \rput(-0.5,2){$2n$}
  \rput(-0.5,3){$3n$}
  \rput(-0.5,4){$4n$}
  \rput(1,-0.4){$1$}
    \rput(2,-0.4){$2$}
    \rput(3,-0.4){$3$}
   \rput(4,-0.4){$4$}
    \rput(5,-0.4){$5$}
    \rput(6,-0.4){$6$}
    \rput(7,-0.4){$7$}
    \rput(8,-0.4){$8$}
   \psaxes[linewidth=.7\pslinewidth]{-}(0,0)(0,0)(8,4)
  \psgrid[griddots=5, subgriddiv=0, gridlabels=0pt](0,0)(8,4)
   \psline(0,4)(3,2)(5,1)(8,0)
\end{pspicture}
 \caption{\label{prang01}$p$-rank $0$ first case}
\end{figure}

There are three Newton polygons for which $e=1$. One is represented in Figure \ref{prang01}.

This is the Newton polygon of  $p(t)$ if and only if $v_p(a_1)\geq n/3$, $v_p(a_2)\geq 2n/3$, $v_p(a_3)=n$ and $v_p(a_4)\geq 3n/2$.  If these conditions hold,   $e=1$  if and only if $p(t)$ has no root in $\mathbb{Q}_p$.

\medskip

\begin{figure}[h]
   \centering
     \begin{pspicture}(0,0)(8,4)
  \psset{unit=7mm,labels=none}
    \rput(-0.4,-0.4){$0$}
  \rput(-0.5,1){$n$}
  \rput(-0.5,2){$2n$}
  \rput(-0.5,3){$3n$}
  \rput(-0.5,4){$4n$}
  \rput(1,-0.4){$1$}
    \rput(2,-0.4){$2$}
    \rput(3,-0.4){$3$}
   \rput(4,-0.4){$4$}
    \rput(5,-0.4){$5$}
    \rput(6,-0.4){$6$}
    \rput(7,-0.4){$7$}
    \rput(8,-0.4){$8$}
   \psaxes[linewidth=.7\pslinewidth]{-}(0,0)(0,0)(8,4)
  \psgrid[griddots=5, subgriddiv=0, gridlabels=0pt](0,0)(8,4)
   \psline(0,4)(4,1)(8,0)
\end{pspicture}
 \caption{\label{prang02}$p$-rank $0$ second case}
\end{figure}
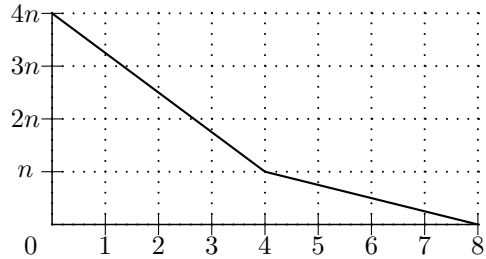

 The second Newton polygon is represented in Figure \ref{prang02}.

This is the Newton polygon of  $p(t)$ if and only if $v_p(a_1)\geq n/4$, $v_p(a_2)\geq n/2$, $v_p(a_3)\geq 3n/4$ and $v_p(a_4)=n$.  If these conditions hold,  $e=1$ if and only if $p(t)$ has no  factor of degrees $1$, $2$ and $3$ in $\mathbb{Q}_p$.

\medskip

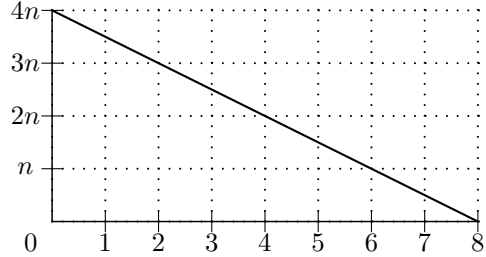
\begin{figure}[h]
   \centering
   \begin{pspicture}(0,0)(8,4)
  \psset{unit=7mm,labels=none}
    \rput(-0.4,-0.4){$0$}
  \rput(-0.5,1){$n$}
  \rput(-0.5,2){$2n$}
  \rput(-0.5,3){$3n$}
  \rput(-0.5,4){$4n$}
  \rput(1,-0.4){$1$}
    \rput(2,-0.4){$2$}
    \rput(3,-0.4){$3$}
   \rput(4,-0.4){$4$}
    \rput(5,-0.4){$5$}
    \rput(6,-0.4){$6$}
    \rput(7,-0.4){$7$}
    \rput(8,-0.4){$8$}
   \psaxes[linewidth=.7\pslinewidth]{-}(0,0)(0,0)(8,4)
  \psgrid[griddots=5, subgriddiv=0, gridlabels=0pt](0,0)(8,4)
   \psline(0,4)(8,0)
\end{pspicture}
 \caption{\label{super}Supersingular case}
\end{figure}

The last Newton polygon is represented in Figure \ref{super}; the corresponding abelian varieties are supersingular.

This is the Newton polygon of  $p(t)$ if and only if $v_p(a_1)\geq n/2$, $v_p(a_2)\geq n$, $v_p(a_3)\geq 3n/2$ and $v_p(a_3)\geq 2n$.  If these conditions hold, $e=1$ if and only if $p(t)$ has no root nor factor of degree $3$ in $\mathbb{Q}_p$.

\section{Supersingular case}
In \cite{Vijay}, Singh, McGuire and Zaytsev gave the list of irreducible characteristic polynomials of supersingular abelian varieties of dimension $4$ , where $q$ is not a square. Here we derive the list when $q$ is a square. 

Let $p(t)$ be an irreducible supersingular Weil polynomial of degree $8$, where $q$ is a square. By Honda-Tate Theorem,  $\frac{1}{q^4}p(\sqrt{q}t)$ is a cyclotomic polynomial of degree $8$ i.e; $\frac{1}{q^4}p(t\sqrt{q})=\Phi_m(t)$ such that $\phi(m)=8$ or $m \in \{15,16,20,24,30\}$. 
Therefore for  each  $m$ above, $p(t)= q^g\Phi_m(\frac{t}{\sqrt{q}})$ gives a supersingular Weil polynomial of degree $8$.
 Let $p(t)=\displaystyle \prod_{i} p_i(t)$ be the decomposition  in irreducible factors of $p(t)$ over $\mathbb{Q}_p$ with $\pi=\sqrt{q}\zeta_n$  as a root, where $\zeta_n$ is primitive $n$th root of unity. To determine the dimension of the corresponding abelian variety, recall  from \cite{milwat}, $p(t)^e$ is a characteristic polynomial of an abelian variety of dimension $4e$, where $e$ is the least common denominator of $\displaystyle \frac{v_p(\pi)}{v_p(q)} \deg p_i(t)= \frac{\deg p_i}{2}$.
Since   $p(t)= q^g\Phi_m(\frac{t}{\sqrt{q}})$,  $\deg p_i=\deg r_i$ where $\Phi_m(t) = \displaystyle\prod_{i}r_i (t)$. But from chapter IV.4 in \cite{Serre}, we have $\deg r_i= r $ where $r$ is the multiplicative order of $p$ in $(\frac{\mathbb{Z}}{m\mathbb{Z}})^*$.\\
Hence, $e=1$ if $r$ is even. In each case of $m$ above, since $\phi(m)=2^3$, $r$ is either even or $r=1$. The later case only happens when $p \equiv 1 \mod m$. 
Hence, $p(t)= q^g\Phi_m(\frac{t}{\sqrt{q}})$, where $p \not\equiv 1 \mod m $ is an irreducible characteristic polynomial of a supersingular abelian variety of dimension 4, for each  $m \in \{15,16,20,24,30\}$.

\end{document}